\documentclass[brochure,english,12pt]{bourbaki}
\usepackage[matrix,arrow]{xy}
\usepackage{amssymb,amsfonts,amsmath,footnote}
\usepackage[francais]{babel}
\usepackage{tikz}
\usetikzlibrary{arrows}
\usepackage{euscript}

    \def\HC{{\mathcal{H}}}

    \def\LC{{\mathcal{L}}}

    \def\XC{{\mathcal{X}}}
    \def\YC{{\mathcal{Y}}}

\def\FS{{\EuScript F}}

\def\LS{{\EuScript L}}

\def\a{\alpha}
\def\b{\beta}

\def\e{\varepsilon}

\def\l{\lambda}
\def\L{\Lambda}
\def\O{\Omega}

\def\z{\zeta}

\newcommand{\into}{\hookrightarrow}
\newcommand{\onto}{\twoheadrightarrow}
\newcommand{\simto}{\stackrel{\sim}{\to}}
\newcommand{\triright}{\stackrel{[1]}{\to}}
\newcommand{\ptau}{{}^p\tau}

\newcommand{\pH}{{}^p \mathcal{H}}

\newcommand{\pD}{{}^p D}

\newcommand{\QM}{\mathbb{Q}}
\newcommand{\RM}{\mathbb{R}}
\newcommand{\CM}{\mathbb{C}}
\newcommand{\ZM}{\mathbb{Z}}
\newcommand{\PM}{\mathbb{P}}
\newcommand{\DM}{\mathbb{D}}

\DeclareMathOperator{\Hom}{Hom}
\DeclareMathOperator{\End}{End}
\DeclareMathOperator{\gr}{gr}
\DeclareMathOperator{\im}{im}
\DeclareMathOperator{\rad}{rad}
\DeclareMathOperator{\supp}{supp}

\date{Mars 2016}
\bbkannee{68\`eme ann\'ee, 2015-2016}
\bbknumero{1115}
\title{The Hodge theory of the Decomposition Theorem}
\subtitle{after M. A. de Cataldo and L. Migliorini}
\author{Geordie WILLIAMSON}
\address{Max--Planck--Institut f\"ur Mathematik \\
Vivatgasse 7\\
D--53111 Bonn, Germany}
\email{geordie@mpim-bonn.mpg.de }

\begin{document}
\maketitle

\noindent{\bf INTRODUCTION}

The Decomposition Theorem is a beautiful theorem about algebraic maps. In
the words of MacPherson \cite{Mac}, \og it contains
as special cases the deepest homological properties of algebraic maps
that we know.\fg{}
Since its proof in 1981 it has found spectacular applications in number
theory, representation theory and combinatorics. Like its cousin the
Hard Lefschetz Theorem, proofs appealing to the Decomposition Theorem are
usually difficult to obtain via other means. This leads one to regard
the Decomposition Theorem as a deep statement lying at the 
heart of diverse problems.

The Decomposition Theorem was first proved by Beilinson, Bernstein,
Deligne and Gabber \cite{BBD}. Their proof proceeds by reduction to positive
characteristic in order to use the Frobenius endomorphism and its
weights, and ultimately rests on Deligne's proof of the Weil
conjectures. Some years later Saito obtained another proof of the
Decomposition Theorem as a corollary of his theory of mixed Hodge
modules \cite{S1,S2}. Again the key is a notion of weight.

More recently, de Cataldo and Migliorini discovered a simpler
proof of the Decomposition Theorem \cite{dCM,dCM2}. The proof is an ingenious
reduction to statements about the
cohomology of smooth projective varieties, which they establish via
Hodge theory. In their
proof they uncover several remarkable geometric statements which go a long way
to explaining \og why\fg{} the Decomposition Theorem holds, purely in the
context of the topology of algebraic varieties. For
example, their approach proves that the intersection
cohomology of a projective variety is of a motivic nature (\og Andr\'e
motivated\fg{}) \cite{dCM5}. Their techniques were
adapted by Elias and the author to prove the existence of Hodge
theories attached to Coxeter systems (\og Soergel modules\fg{}),
thus proving the Kazhdan-Lusztig positivity conjecture \cite{EW}.

The goal of this article is to provide an overview of the main ideas
involved in de Cataldo and
Migliorini's proof. A striking aspect of the proof is that
it gathers the Decomposition Theorem together with several
other statements generalising the Hard Lefschetz Theorem and the
Hodge-Riemann Bilinear Relations (the \og
Decomposition Theorem Package\fg{}). Each ingredient is indispensable
in the induction. One is left with the impression that the
Decomposition Theorem is not a theorem by itself, but rather belongs
to a family of statements, each of which sustains the others.



Before stating the Decomposition Theorem we recall two earlier
theorems concerning the topology of algebraic maps. The first
(Deligne's Degeneration Theorem) is an instance of the Decomposition Theorem. 
The second (Grauert's Theorem) provides an illustration of the appearance of a
definite form, which eventually forms part of the \og Decomposition Theorem
Package\fg{}.

\subsection{Deligne's Degeneration Theorem} Let $f : X \to Y$ be a smooth
(i.e. submersive) 
projective morphism of complex algebraic varieties. Deligne's theorem
asserts that the Leray spectral sequence
\begin{equation}
E^{pq}_2 = H^p(Y, R^q f_* \QM_X) \Rightarrow H^{p + q}(X,
\QM)\label{eq:E2}
\end{equation}
is degenerate (i.e. $E_2 = E_\infty$). Of course such a statement is
false for submersions between manifolds (e.g. the Hopf fibration).
The theorem asserts that something very special happens for smooth
algebraic maps.

Let us recall how one may construct the Leray spectral sequence. In
order to compute the cohomology of $X$ we replace the constant
sheaf $\QM_X$ on $X$ by an injective resolution. Its direct image on
$Y$ then has a natural \og truncation\fg{} filtration whose successive
subquotients are the (shifted) higher direct image sheaves $R^qf_*
\QM_X[-q]$. This filtered complex of sheaves gives rise to the Leray
spectral sequence.

In fact, Deligne proved that there exists a decomposition in
the derived category of sheaves on $Y$:
\begin{equation}
Rf_* \QM_X \cong \bigoplus_{q \ge 0} R^qf_* \QM_X[-q]\label{eq:2}
\end{equation}
 (i.e. the filtration of the previous paragraph splits).
The decomposition in \eqref{eq:2} implies the degeneration of
\eqref{eq:E2}, and in fact is the universal explanation for such a
degeneration. Deligne also proved that each 
local system $R^qf_* \QM_X$ is semi-simple. Hence the object $Rf_*
\QM_X$ is as semi-simple as we could possibly hope. This is the
essence of the Decomposition Theorem, as we will see.

Because $f : X \to Y$ is smooth and projective any fibre of $f$ is
a smooth projective variety. Deligne deduces the decomposition in
\eqref{eq:2} by applying the Hard Lefschetz Theorem to the cohomology of the
fibres of
$f$. Thus the decomposition of $Rf_*\QM_X$ is deduced
from a deep fact about the global cohomology of a smooth projective
variety. This idea occurs repeatedly in the proof of de Cataldo
and Migliorini.

\subsection{Grauert's Theorem}

Let $X$ denote a smooth projective surface and let $C = \bigcup_{i =
  1}^m C_i$ denote a connected union of irreducible curves on $X$. It
is natural to ask whether $C$ can be contracted. That is, whether there exists a map
\[
f : X \to Y
\]
which is an isomorphism on the complement of $C$ and contracts $C$ to a
point. Of course such a map of topological spaces always exists, but
it is a subtle question if one requires $f$ and $Y$ to be algebraic or
analytic. An answer is given by Grauert's theorem: $f$ exists
analytically if and only if the intersection form
\begin{equation}
( [C_i] \cap [C_j] )_{1 \le i, j \le k}\label{eq:gnd}
\end{equation}
is negative definite. For example, if $C$ is irreducible (i.e. $k=1$)  
then $C$ can be contracted if and only if $C$ has
negative self-intersection in $X$.

Let us assume that such an $f$ exists, and let $y \in Y$ denote the
image of $C$. Then in this case the
Decomposition Theorem asserts a decomposition in the derived category
of sheaves on $Y$
\begin{equation}
  \label{eq:4}
  Rf_* \QM_X[2] = IC(Y) \oplus \bigoplus_{i = 1}^k \QM_{y}
\end{equation}
where $IC(Y)$ is a complex of sheaves on $Y$ which is a simple object in
the category of perverse sheaves. Again, \eqref{eq:4} can be
interpreted in the language of perverse
sheaves as saying that the object $Rf_* \QM_X[2]$ is as semi-simple
as possible.

Remarkably, the decomposition in \eqref{eq:4} is equivalent to
the fact that the intersection form in \eqref{eq:gnd} is
non-degenerate. Thus in this example the Decomposition Theorem is a
consequence of a topological fact about contractibility of
curves on a surface.
Note also that here the geometric theorem that we are
using (negative definiteness) is stronger
than what we need for the Decomposition Theorem
(non-degeneracy). As we will see, keeping track of such signs
plays an important role in de Cataldo and Migliorini's proof.

\subsection{Structure of the Paper}

This paper consists of three sections. In \S 1 we recall
the necessary background from topology, Hodge theory and perverse sheaf theory
and state the Decomposition Theorem. In \S 2 we discuss de Cataldo and
Migliorini's proof for semi-small maps. The case of semi-small maps
has the advantage of illustrating several of the general features of
the proof very well,  whilst being much simpler in structure. In \S 3
we give the statements and an outline of the main steps of the induction 
establishing the theorem for arbitrary maps.

\bigskip

\emph{Acknowledgements} --- I would like to thank M. A. de
Cataldo, D. Juteau,
C. Mautner, L. Migliorini, W. Soergel and K. Vilonen for
many interesting conversations about perverse sheaves. In addition,
thanks to F. El Zein, S. Riche, G. Sacc\`a and J. Torres for comments on a
preliminary version.

\section{Background} \label{sec:background}

In this section we briefly recall the tools (intersection forms, classical Hodge theory,
perverse sheaves) which we will be using
throughout this paper. We discuss the relationship between
perverse sheaves and the weak Lefschetz theorem and state the
Decomposition Theorem.

\begin{rema}
  A remark on coefficients: The natural setting for the Decomposition Theorem and
  its relatives is that of sheaves of $\QM$-vector spaces. However, at
  some points below it is necessary to consider sheaves of $\RM$-vector spaces
  (usually due to limit arguments). To avoid repeated change of
  coefficients we have chosen to work with $\RM$-coefficients
  throughout. All of the arguments of this paper are easily adapted
  for $\QM$-coefficients,  as the reader
  may readily check.
\end{rema}

\subsection{Algebraic Topology} \label{sec:at}

All spaces will be complex algebraic varieties equipped with their
classical (metric) topology. The dimension of a complex
algebraic variety will always mean its complex dimension. We do not
assume that varieties are irreducible, and dimension means
the supremum over the dimension of its components. Given a variety $Z$ we denote by
\[
H^*(Z)  = H^*(Z,\RM) \quad \text{and} \quad  H_*(Z) = H_*(Z,\RM)
\]
its singular cohomology and singular homology with closed supports (``Borel-Moore
homology''), with coefficients in the real numbers.

Any irreducible subvariety $Z' \subset Z$ of dimension $p$ has a
fundamental class
\[
[Z'] \in H_{2p}(Z).
\]
If $Z$ is of dimension $n$ then $H_{2n}(Z)$ has a basis given by the
fundamental classes of irreducible components of maximal dimension.

If $X$ is smooth of dimension $n$ then (after choosing once and for all
an orientation of $\CM$) Poincar\'e duality
gives a canonical isomorphism
\begin{equation}
H_p(X) \simto H^{2n-p}(X).\label{eq:PD}
\end{equation}
If $X$ is in addition compact then $H^*(X)$ has a
non-degenerate \emph{Poincar\'e form}
\[ ( -, -) : H^{2n-p}(X) \times H^p(X) \to \RM \]
and $H_*(X)$ is equipped with a non-degenerate \emph{intersection
  form}
\[
\cap : H_p(X) \times H_q(X) \to H_{p+q-2n}(X).
\]
These forms match under Poincar\'e duality. If $X$ is smooth we will often identify 
$H^*(X)$ with the real de Rham cohomology of $X$. 
In de Rham cohomology the Poincar\'e form is given by the integral
\[
(\a, \b) \mapsto \int_X \a \wedge \b.
\]

Suppose $Z$ is a proper closed subvariety inside a smooth
$n$-dimensional variety $X$. If $p + q = 2n$ the
inclusion $Z \into X$ gives rise to an intersection form (see e.g.
\cite[Chapter 19]{F})
\[
H_p(Z) \times H_q(Z) \to \RM.
\]
Geometrically this corresponds to moving cycles on $Z$ into $X$ until
they become transverse, and then counting the number of intersection
points. If $X$ is proper and connected the map $H_*(Z) \to H_*(X)$ is an
isometry for intersection forms.

\subsection{Hodge Theory}

Let $X$ be a smooth and connected projective variety of complex dimension
$n$. Let $H^*(X)$ denote the de Rham cohomology of $X$ with
coefficients in the real numbers. Throughout it will be convenient to
shift indices;  consider the finite-dimensional
graded vector space
\[
H = \bigoplus_{i \in \ZM} H^i \quad \text{where} \quad H^i := H^{n+i}(X).
\]
Under this normalization the Poincar\'e pairing induces canonical
isomorphisms
\begin{equation}
  \label{eq:1}
  H^{-i} \simto (H^i)^\vee \quad \text{for all $i \in \ZM$}
\end{equation}
where $(H^i)^\vee$ denotes the dual vector space.

\begin{theo} [The Hard Lefschetz Theorem]
Let $\omega \in H^2(X)$ denote the
  Chern class of an ample line bundle. For all $i \ge 0$,
  multiplication by $\omega^i$ induces an isomorphism
  \begin{equation}
    \label{eq:3}
    \omega^i : H^{-i} \simto H^i.
  \end{equation}
\end{theo}

Let $P^{-i} \subset H^{-i}$ denote the primitive subspace:
\[
P^{-i} := \ker ( \omega^{i+1} : H^{-i} \to H^{i+2}).
\]
The Hard Lefschetz Theorem gives the primitive decomposition:
\[
\bigoplus_{i \ge 0} \RM[\omega]/(\omega^{i+1}) \otimes_{\RM} P^{-i}
\simto H.
\]

\begin{rema} \label{rem:sl2}
  Consider the Lie algebra $\mathfrak{sl}_2 := \RM f \oplus \RM h
  \oplus \RM e$ with
\[
f = \left ( \begin{matrix} 0 & 0 \\ 1 &
      0 \end{matrix} \right), \quad
h = \left ( \begin{matrix} 1 & 0 \\ 0 & -1
       \end{matrix} \right) \quad \text{and} \quad 
e = \left ( \begin{matrix} 0 & 1 \\ 0 &
      0 \end{matrix} \right).
\]
The Hard Lefschetz Theorem is equivalent to the existence of a
$\mathfrak{sl}_2$-action on $H$ with $e(x) = \omega \wedge x$ and $h(x) = j x$ for
all $x \in H^j$. The primitive decomposition is the isotypic
decomposition and the primitive subspaces
are the lowest weight spaces.
\end{rema}

We now state the Hodge-Riemann bilinear relations, for which we need a little more
notation. For $i \ge 0$ the form
\[
Q(\a, \b) := \int \omega^i \wedge \a \wedge \b
\]
on $H^{-i}$ is symmetric if $n-i$ is even and alternating if $n-i$ is
odd. It is non-degenerate by the Hard Lefschetz theorem. Given a real
vector space $V$ we denote by $V_\CM$ its 
complexification. The form
\[
\kappa(\a, \b) := (\sqrt{-1})^{n-i} Q(\a,\overline{\b})
\]
on $H^{-i}_\CM$ is Hermitian and non-degenerate.

Consider the Hodge decomposition and corresponding
primitive spaces
\[
H^j_{\CM} = \bigoplus_{p + q = n + j} H^{p,q}, \qquad P^{p,q} := P^{p+q-n}_{\CM} \cap H^{p,q}.
\]

\begin{theo}[Hodge-Riemann Bilinear Relations] \label{thm:HR}
The Hodge decomposition
  is orthogonal with respect to $\kappa$. Moreover, if  $\a \in
  P^{p,q}$ is non-zero and $k := p + q$ then
\[
(\sqrt{-1})^{p-q-k}(-1)^{k(k-1)/2}\kappa(\a,\overline{\a}) > 0.
\]
\end{theo}

\begin{rema} The Hodge-Riemann
  relations imply that the restriction of the Hermitian form
  $\kappa$ to $P^{p,q}$ is definite of a fixed sign. This fact is
  crucial below. As long as the reader keeps this definiteness in mind, the
  precise nature of the signs can be ignored on a first reading.
\end{rema}

\begin{rema} \label{rem:amplecone}
  More generally, hard Lefschetz and the Hodge-Riemann bilinear
  relations are valid for any class $\omega \in H^2(X)$ in
  the ample cone (the convex hull of all strictly positive real multiples
  of ample classes).
\end{rema}

A \emph{(real, pure) Hodge structure} of weight $k$ is a
finite-dimensional real vector
space $V$ together with a decomposition $V_{\CM} =
\bigoplus_{p + q = k} V^{p,q}$ such that $\overline{V^{p,q}} =
V^{q,p}$. Hodge structures form an abelian category in a natural
way. A \emph{polarisation} of a real Hodge structure of weight $k$ is
a bilinear form $Q$ on $V$ which is symmetric if $k$ is even,
anti-symmetric if $k$ is odd and such that the corresponding
Hermitian form $\kappa(\a, \b) := (\sqrt{-1})^k Q(\a, \overline{\b})$
on $V_{\CM}$ satisfies the Hodge-Riemann Bilinear Relations (Theorem
\ref{thm:HR}). For example, for $i \ge 0$ each $H^{-i}$ above is a
Hodge structure of weight $n-i$ and $P^{-i} \subset H^{-i}$ is a
Hodge substructure polarised by $Q$.

\subsection{Constructible and Perverse Sheaves}

In the following we recall the formalism of the constructible derived
category. For more detail the reader is referred to \cite[\S 5]{dCM4} and the references
therein.

We denote by $D^b_c(Y)$ the constructible derived category of sheaves
of $\RM$-vector spaces on $Y$. This is a triangulated category with
shift functor $[1]$.
Given an object $\FS \in D^b_c(Y)$ we
denote by $\HC^i(\FS)$ its cohomology sheaves.  Given a morphism $f :
X \to Y$ of algebraic varieties we have functors
\begin{equation*}
  \begin{tikzpicture}[xscale=2]
    \node (x) at (-1,0) {$D^b_c(X)$};
    \node (y) at (1,0) {$D^b_c(Y)$};
\draw[->] (x) .. controls (0,0.6) .. node[above] {$f_*, f_!$} (y);
\draw[<-] (x) .. controls (0,-0.6) .. node[below] {$f^*, f^!$} (y); 
 \end{tikzpicture}
\end{equation*}
(we only consider derived functors and write $f_*$ instead of $Rf_*$
etc.). Verdier duality is denoted $\DM : D^b_c(Y) \to
D^b_c(Y)$. 

We let $\RM_Z$ and $\omega_Z$ denote the constant and dualizing
sheaves on $Z$. If $Y$ is smooth we have $\omega_Y
= \RM_Y[2\dim Y]$ canonically (Poincar\'e duality). Given $\FS \in D^b_c(Y)$ we denote
its hypercohomology by $H(Y,\FS)$. In the notation of \S~\ref{sec:at} we have
\[
H^j(Y) = H^j(Y, \RM_Y) \quad \text{and} \quad H_j(Y) = H^{-j}(Y, \omega_Y).
\]

The full subcategories
\begin{gather*}
  \pD^{\le 0}(Y) := \{ \FS \in D^b_c(Y) \; | \;
  \dim \supp \HC^i(\FS) \le -i \text{ for all $i$}\}, \\
  \pD^{\ge 0}(Y) :=   \{ \FS \in D^b_c(Y) \; | \;
  \dim \supp \HC^i(\DM\FS) \le -i \text{ for all $i$}\}
\end{gather*}
define a t-structure on $D^b_c(Y)$ whose heart is the abelian category
$P_Y \subset D^b_c(Y)$ of perverse sheaves (for the middle perversity). (The standard warning that perverse sheaves
are not sheaves, but rather complexes of sheaves is repeated here.)

Define $\pD^{\le m} := \pD^{\le 0}[-m]$ and $\pD^{\ge m} :=
  \pD^{\ge 0}[-m]$. We denote by $\ptau_{\le m}$ and $\ptau_{\ge m}$ the perverse
truncation functors
\[
\ptau_{\le m} : D^b_c(Y) \to \pD^{\le m}(Y) \quad \text{ and } \quad
\ptau_{\ge m} : D^b_c(Y) \to \pD^{\ge m}(Y)
\]
which are right (resp. left) adjoint to the inclusion functors. Given
$\FS \in D^b_c(Y)$ its perverse cohomology groups are $\pH^i(\FS) :=
\ptau_{\le 0} \ptau_{\ge 0}(\FS[i]) \in P_Y$.

\begin{rema}
  For fixed $\FS \in D^b_c(Y)$ we have $\ptau_{\le i} \FS = 0$ for $i
  \ll 0$ and $\ptau_{\ge i} \FS = 0$ for $i \gg 0$. It is 
  convenient to view $\FS$ as equipped with a canonical
  exhaustive filtration (in the triangulated sense)
\[
\dots \to \ptau_{\le i} \FS \to \ptau_{\le i+1} \FS \to \dots
\]
with subquotients the (shifted) perverse sheaves $\pH^i(\FS)[-i]$.
\end{rema}

Given any locally closed, smooth and connected subvariety $Z \subset
Y$ and a local system $\LS$ of $\RM$-vector spaces on $Z$ we
denote by $IC(\overline{Z}, \LS)$ the intersection
cohomology complex of $\LS$. The object $IC(\overline{Z}, \LS) \in P_Y$ is
simple if $\LS$ is, and all simple perverse sheaves are of this form. For example,
if $\overline{Z}$ is smooth and $\LS$ extends as a local system
$\overline{\LS}$ to $\overline{Z}$ then $IC(\overline{Z}, \LS) =
\overline{\LS}[\dim Z]$.
We write $IH(\overline{Z}, \LS) = H(Y, IC(\overline{Z}, \LS))$ for the
intersection cohomology of $\overline{Z}$ with coefficients in
$\LS$. 
 If $\LS$ is the trivial local system we write
$IC(\overline{Z})$ and $IH(\overline{Z})$ instead of $IC(\overline{Z},
\LS)$ and $IH(\overline{Z}, \LS)$.

Let us fix a Whitney stratification $Y = \bigsqcup_{\l \in \Lambda}
Y_\l$ and denote by $i_\mu : Y_\mu \into Y$ the inclusion. If we fix a
stratum $Y_\l \subset Y$ and a local system
$\LS$ on $Y_\l$ then $IC(\overline{Y_\l}, \LS)$ is uniquely characterised by the conditions:
\begin{gather}
i_\l^* IC(\overline{Y}_{\lambda}, \LS) = \LS[\dim Z], \\
\HC^j(i_\mu^* IC(\overline{Y}_\l, \LS)) = 0 \quad \text{for $j \ge
  -\dim Y_\mu$ and $\mu \ne
  \lambda$,}\\
\HC^j(i_\mu^! IC(\overline{Y}_\l, \LS)) = 0 \quad \text{for $j \le
  -\dim Y_\mu$ and $\mu \ne \lambda$}.
\end{gather}

At several points in de Cataldo and Migliorini's proof vanishing
theorems for perverse sheaves on affine varieties play an important
role. Recall Artin-Grothendieck vanishing (see e.g. \cite[3.1.13]{La}): if $\FS$ is a
constructible \emph{sheaf} (i.e. $\FS = \HC^0(\FS)$) on an affine
variety $U$ then
\begin{equation}
  \label{eq:affinevanish}
  H^j(U, \FS) = 0 \quad \text{for $j > \dim U$.}
\end{equation}
The following proposition characterises the perverse sheaves as those
complexes for which such vanishing is universal:

\begin{prop} \label{prop:affinevanish}
$\FS \in D^b_c(Y)$ belongs to $\pD^{\le 0}(Y)$ if and only if, for all
affine open subvarieties $U \subset Y$, we have
\[
H^j(U, \FS) = 0 \quad \text{for $j > 0$}.
\]
Similarly, $\FS \in \pD^{\ge 0}$ if and only if for all affine open $U$ we have
\[
H^j_c(U, \FS) = 0 \quad \text{for $j < 0$},
\]
where $H^j_c(U, \FS)$ denotes cohomology with compact
supports.
\end{prop}

\noindent{\sc Proof} (Sketch) --- The first statement implies the
second, by Verdier duality. The implication $\Rightarrow$ is
easily deduced from the definition of the perverse $t$-structure
and Artin-Grothendieck vanishing \eqref{eq:affinevanish}. For the 
implication $\Leftarrow$ see \cite[4.1.6]{BBD}.

\bigskip

Now suppose that $Y$ is projective and let $i : D \into Y$ denote the
inclusion of a hyperplane section and $j : Y \setminus D \into Y$ the open
inclusion of its (affine) complement. After taking cohomology of the
distinguished triangle $j_!j^!\FS \to \FS \to i_*i^*\FS \triright$ 
or its dual and applying the above vanishing we deduce:

\begin{theo} [Weak Lefschetz for Perverse Sheaves] \label{prop:perversewl}
Let $\FS \in D^b_c(Y)$ be perverse.
  \begin{itemize}
  \item The restriction map $H^j(Y, \FS) \to H^j(D, i^*\FS)$ is an
    isomorphism for $j < -1$ and is injective for $j = -1$.
  \item The pushforward map $H^j(D, i^!\FS) \to H^j(Y, \FS)$ is an
    isomorphism for $j > 1$ and is surjective for $j = 1$.
  \end{itemize}
\end{theo}

\subsection{The Decomposition Theorem} 

\begin{defi} 
An object in $D^b_c(Y)$ is \emph{semi-simple} if it is isomorphic to a direct
sum of shifts of intersection cohomology complexes of semi-simple
local systems.
\end{defi}

\begin{theo}[Decomposition Theorem]
If $f : X \to Y$ is projective and $X$ is smooth then $f_*\RM_X$ is semi-simple.
\end{theo}

\begin{rema} Some remarks concerning the generality of the
  Decomposition Theorem discussed below:
  \begin{itemize}
\item One could drop the assumption that $X$ be smooth and
  replace $f_*\RM_X$ by $f_* IC(X)$. This formulation follows from the above
 via resolution of singularities. By Chow's lemma
  we could also replace \og $f$ projective\fg{} by \og $f$
  proper\fg{}. The formulation above is preferred because it is the
  one addressed in this paper.
  \item In Saito's theory the Decomposition Theorem is proved more
    generally for $f_* IC(X, \LS)$ where $\LS$ is any local system
    underlying a polarisable variation of Hodge structure on a Zariski
    open subvariety of $X$. It is likely that the techniques
    discussed here could handle this case (after reducing to the
    normal crossing situation and using
    the existence of a pure Hodge structure on $IH(X,\LS)$
    established by Kashiwara-Kawai \cite{KK}, and
    Cattani-Kaplan-Schmid \cite{CKS}). Recently El Zein, L\^e and Ye have 
    proposed another proof of the Decomposition Theorem in this level of
    generality \cite{EZ1,EZ2,EZ3}.
 \item More general still are the results of Sabbah \cite{Sab} and
    Mochizuki \cite{M}
    which establish the semi-simplicity of $f_*IC(X, \LS)$ where
    $\LS$ is \emph{any} semi-simple  $\CM$-local system. The proof is
    via a generalization of Saito's theory, and probably goes far
    beyond what is possible with the techniques discussed here.
  \end{itemize}
\end{rema}


\bigskip
\section{Semi-small maps} \label{sec:ss}

\subsection{The Decomposition Theorem for Semi-Small Maps}

Suppose (as we will assume throughout this paper) that $X$ is smooth, connected and projective of
complex dimension $n$ and that $f : X
\to Y$ is a surjective algebraic map. Throughout we fix a stratification
\[
Y = \bigsqcup_\Lambda Y_\lambda
\]
of $Y$ adapted to $f$. In particular, each $Y_\lambda$ is connected
and, over each stratum, $f : f^{-1}(Y_\lambda) \to Y_\l$ is a topologically locally trivial
fibration in (typically singular) varieties.


\begin{defi} The map $f$ is \emph{semi-small} if for all $\lambda \in
  \Lambda$ and some (equivalently all) $y \in
  Y_\lambda$ we have
\begin{equation} \label{ss:cond}
\dim f^{-1}(y) \le \frac{1}{2} ( \dim Y - \dim Y_\lambda).
\end{equation}
\end{defi}

Semi-small maps play an important role in the theory of perverse
sheaves. This is mainly because of the following fact (which is a
straightforward consequence of
the proper base change theorem and the Verdier self-duality of $f_*\RM_X[n]$):

\begin{prop} If $f$ is semi-small then $f_*\RM_X[n]$ is perverse.
\label{prop:ss}
\end{prop}

\begin{rema} From the definition it follows that a semi-small
  map is finite on any open stratum of $Y$. It can be useful to think of semi-small
  maps as being the finite maps of perverse  sheaf theory. (Compare
  with the fact that the (derived) direct image
  of the constant sheaf along a projective morphism is a sheaf if and
  only if $f$ is finite.)
\end{rema}

\begin{theo}[Decomposition Theorem for Semi-Small Maps]  \label{thm:dtss}
If $f$ is semi-small then $f_*\RM_X[n]$ is a semi-simple perverse
sheaf. More precisely:
\begin{itemize}
\item We have a
  canonical decomposition
  \begin{equation} \label{eq:dtss}
f_* \RM_X[n] = \bigoplus_{\lambda \in \Lambda} IC(\overline{Y}_\lambda,
\LS_\lambda)
\end{equation}
where each $\LS_\lambda$ is the local system on $Y_\lambda$ associated to
$y \mapsto H^{\dim Y - \dim Y_\lambda}(f^{-1}(y))$.
\item Each local system $\LS_\lambda$ is semi-simple.
\end{itemize}
\end{theo}

\begin{rema}
The semi-small case is special because the decomposition
\eqref{eq:dtss} is canonical and explicit. For general maps the
decomposition is not canonical and it
is difficult to say a priori which summands occur in the
direct image.
\end{rema}

\begin{rema}
An important aspect of the Decomposition Theorem (already 
non-trivial in Deligne's Degeneration Theorem) is that each local system $\LS_\l$
is semi-simple. In the semi-small case the representations
corresponding to each $\LS_\l$ are dual to the permutation
representation of $\pi_1(Y_\l,y)$ on the irreducible components of
the fibre $f^{-1}(y)$ of \og maximal\fg{} (i.e $=\frac{1}{2}( \dim Y
- \dim Y_\lambda)$) dimension. In particular, each representation factors over a
finite group and semi-simplicity follows from Maschke's Theorem in finite
group theory.
\end{rema}

\begin{rema}
The decomposition \eqref{eq:dtss} implies 
  that the cohomology of the fibres of $f$ is completely determined by
  the local systems $\LS_\l$ and the singularities of $Y$. Thus much of the topology of $f$ is
  determined by the irreducible components of each fibre, and the
  monodromy along each stratum. This gives a hint as to the
  nature of the Decomposition Theorem.
\end{rema}

\begin{rema}
  The decomposition \eqref{eq:dtss} gives
  a canonical decomposition of cohomology:
  \begin{equation} \label{eq:dtsscoh}
    H^{*+n}(X) = \bigoplus_{\l \in \Lambda} IH^*(\overline{Y}_\l, \LS_\l).
  \end{equation}
In \cite{dCM6} it is shown that this decomposition is motivic
(i.e. given by algebraic cycles in $X \times_Y X$). For example if
$X$ is proper this gives a canonical decomposition of the Chow motive
of $X$ \cite[Theorem 2.4.1]{dCM6}.
\end{rema}

We say that $\omega \in H^2(X)$ is a \emph{semi-small class} if
$\omega$ is the first Chern class of a line bundle $\LC$, some
positive power of which is globally generated and whose global
sections yield a semi-small map $X \to Y$.

\begin{theo}[Hard Lefschetz and Hodge-Riemann for Semi-Small Classes]
 \label{thm:lef}
Let $\omega\in H^2(X)$ be a semi-small class. Then multiplication by $\omega$
  satisfies hard Lefschetz and the Hodge-Riemann bilinear relations.
 \end{theo}

\begin{rema}
More generally one can show that if $f : X \to Y$ is any morphism, $\LC$ is
an ample line bundle on $Y$ and $\omega$ is the Chern class of $f^* \LC$
then $\omega$ satisfies hard Lefschetz if and only if $f$ is
semi-small, see \cite[Proposition 2.2.7]{dCM}.
\end{rema}

\begin{rema}
 If one knows that the hypercohomology of each summand appearing in
 the Decomposition Theorem satisfies hard Lefschetz and
  the Hodge-Riemann relations (as follows for example from Saito's theory) then Theorem
  \ref{thm:lef} is an immediate consequence of Theorem
  \ref{thm:dtss}. A key insight of de Cataldo and Migliorini
  is to realise that the Decomposition Theorem in the semi-small case is implied by
  Theorem \ref{thm:lef}, as we will explain below.
\end{rema}

\begin{rema}
  Theorem \ref{thm:lef} can be used to put pure Hodge structures on
  each summand in \eqref{eq:dtsscoh}.
\end{rema}
\subsection{Local Study of the Decomposition Theorem: Semi-Small Case}

Suppose that $f : X \to Y$ is as in the previous section with $f$
semi-small. From the definition of a semi-small map it is immediate
that the dimension of any fibre of $f$ is at most half of the
dimension of $X$, and that equality can only occur at finitely many points
in $Y$. It is useful to think of these points as the \og most singular
points\fg{} of $f$.

\begin{exem} The first interesting example of a semi-small map is
  that of a contraction of curves on a surface appearing in
  Grauert's theorem (as discussed in the introduction). The image of the contracted curves is
  typically a singular point of $Y$, which is an example of a $y \in
  Y$ that we study below.
\end{exem}

Let us assume that $X$ is of even dimension $n = 2m$. We fix a point
$y \in Y$ such that $\dim f^{-1}(y) =m$. Consider the Cartesian
diagram:
\begin{equation}
  \label{eq:6}
  \xymatrix{ F \ar[r]^i \ar[d]^f & X \ar[d]^f \\
\{ y \} \ar[r]^i & Y}
\end{equation}
The fibre $F = f^{-1}(y)$ is typically reducible. If we denote by $F_1, F_2, \dots, F_k$ the irreducible components of
$F$ of dimension $m$ then we have
\begin{equation}
  \label{eq:7}
  H_n(F) = \bigoplus_{i = 1}^k \RM [F_i]
\end{equation}
where $[F_i] \in H_n(F)$ denotes the fundamental class of $F_i \subset
F$. Because each $F_i$ is half-dimensional inside $X$ the inclusion $F
\into X$ equips $H_n(F)$ with a symmetric intersection form (see \S~\ref{sec:at})
\begin{equation} \label{eq:lif}
H_n(F) \times H_n(F) \to \RM.
\end{equation}
We will call this form the \emph{local intersection form} (at $y$).

The Decomposition Theorem predicts
\begin{equation}
f_* \RM_X[n] = \FS \oplus i_* (H^n(F))_y\label{eq:8}
\end{equation}
where $i_* (H^n(F))_y$ denotes the constant sheaf on $\{ y \}$ with
stalk $H^n(F) = H_n(F)^\vee$. (Here $\FS$ is some perverse sheaf,
whose structure can be ignored for the moment.) We will say that the
\emph{Decomposition Theorem holds at $y$} if the decomposition
\eqref{eq:8} is valid.

\begin{rema} Let us justify this terminology. In de Cataldo and Migliorini's proof one
  knows by induction that the restriction of $\FS$ to the complement
  of all of the point strata is semi-simple. It is then not difficult
  to prove that the decomposition \eqref{eq:8} for all point strata (or \og most
  singular points\fg{}) is equivalent to the Decomposition Theorem
for $f_* \RM_X[n]$. Thus the innocent looking \eqref{eq:8} is the key to the
Decomposition Theorem for Semi-Small Maps.
\end{rema}

How do we decide whether the Decomposition Theorem holds at $y$? The
Decomposition Theorem holds at $y$ if and only if the skyscraper sheaf $i_* \RM_y$ occurs
with multiplicity equal to the dimension of $H_n(F)$ as a summand of
$f_* \RM_X[n]$. We can rephrase this as follows: If we consider the pairing
  \begin{equation}
\Hom(i_* \RM_{y}, f_* \RM_X[n]) \times \Hom(f_* \RM_X[n],i_* \RM_{y})
\to \End(i_* \RM_{y}) = \RM.\label{eq:10}
\end{equation}
then the Decomposition Theorem holds at $y$ if and only if the rank of
the pairing \eqref{eq:10} is $\dim H_n(F)$.

\begin{lemm} \label{lem:caniso}
  We have canonical isomorphisms
  \begin{equation}
    \label{eq:5}
\Hom(i_* \RM_{y}, f_* \RM_X[n]) = H_n(F) = \Hom(f_* \RM_X[n],i_* \RM_{y}).
  \end{equation}
\end{lemm}

\noindent{\sc Proof} --- By adjunction, the proper base
change theorem and the identification $\RM_X[n] = \omega_X[-n]$
(remember that $X$ is smooth) we have 
\begin{gather*}
  \label{eq:9}
  \Hom(i_* \RM_y, f_* \RM[n]) = \Hom(\RM_y, i^!f_* \RM_X[n])
=   \Hom(\RM_y, f_*i^! \RM_X[n]) = \\ = \Hom(\RM_y, f_* \omega_F[-n])
= \Hom(\RM_F, \omega_F[-n])  = H_n(F).
\end{gather*}
The identification $H_n(F) = \Hom(f_* \RM_X[n],i_* \RM_{y})$
follows similarly.
\bigskip

Using the identifications \eqref{eq:5} we can rewrite the form in
\eqref{eq:10} as a pairing
\begin{equation}
H_n(F) \times H_n(F) \to \RM.\label{eq:11}
\end{equation}
The following gives the geometric significance of \eqref{eq:10}
(see \cite{dCM} and \cite[Lemma 3.4]{JMW}):

\begin{lemm} The form \eqref{eq:11} agrees with the local intersection
  form \eqref{eq:lif}.
\end{lemm}

From this discussion we conclude:

\begin{prop} \label{prop:iff}
  The Decomposition Theorem holds at $y$ if and only if the local
  intersection form is non-degenerate.
\end{prop}

\subsection{The Semi-Small Index Theorem}

We keep the notation of the previous section. In particular
\[
f : X \to Y
\]
is a semi-small map, $\dim X = n =  2m$, and $y \in Y$ is such that
$F := f^{-1}(y)$ is of (half) dimension $m$. In the previous section
we outlined a reduction of the Decomposition Theorem in the 
semi-small case to checking that the local intersection form on $H_n(F)$ is
non-degenerate. In fact, a stronger statement is true:

\begin{theo} [Semi-Small Index Theorem, \cite{dCM}]
  The local intersection form on $H_n(F)$ is $(-1)^m$-definite.
\end{theo}

In this section we explain how to deduce this theorem from the
Hodge-Riemann relations for semi-small classes (Theorem \ref{thm:lef}).

\begin{rema} \label{rem:deg0}
  The Semi-Small Index Theorem remains true for any proper semi-small map $f :
  X \to Y$, as long as $X$ is quasi-projective and smooth \cite[Corollary
  2.1.13]{dCM2}. The proof requires the
  Decomposition Theorem (with signs) for an arbitrary map. (One can
  compactify $f : X \to Y$ but one may destroy semi-smallness.)
\end{rema}

\noindent{\sc Proof} (Sketch) --- Consider the composition
\[
cl : H_n(F) \to H_n(X) \simto H^n(X)
\]
where the first map is induced from the inclusion $F \into X$ and the
second map is Poincar\'e duality. The spaces $H_n(F)$ and $H_n(X)$ are
equipped with intersection forms and $H^n(X)$ carries its Poincar\'e
form. By basic algebraic topology:
\begin{equation}
  \label{eq:16}
  cl \text{ is an isometry.}
\end{equation}
We will use the Hodge-Riemann relations for $H^n(X)$ to deduce the
index theorem. The bridge to the Hodge-Riemann relations is provided
by the following two beautiful facts:

\begin{lemm} 
Let $\omega$ denote the Chern class of $f^*\LC$, for $\LC$ an ample line bundle on $Y$. 
The image of $cl$ consists of $\omega$-primitive classes of Hodge type $(m,m)$.
\end{lemm}

\noindent{\sc Proof} --- Recall that $H_n(F)$ has a basis consisting
of fundamental classes $[F_i]$ of irreducible components of $F = f^{-1}(y)$ of
maximal dimension. Thus the image of $cl$ consists of algebraic
cycles, and the claim about Hodge type follows. It remains to see
that the image consists of primitive classes. Under the isomorphism $H_n(X) \simto H^n(X)$
multiplication by $\omega$ on the right corresponds to intersecting
with a general hyperplane section of $f^*\LC$ on the left. We may
assume that such a hyperplane section is the inverse image, under $f$,
of a general hyperplane section of $\LC$. However such a hyperplane
section has empty intersection with $\{ y \}$ (being a point) and
hence its inverse image does not intersect $F$. The claim follows.

\begin{lemm} \label{lem:injective}
  $cl$ is injective.
\end{lemm}

\noindent{\sc Proof} --- The pushforward $H_n(F) \to H_n(X)$ is dual to the
restriction map
\[
r : H^n(X) \to H^n(F).
\]
We will show that $r$ is surjective, which implies the lemma.

Let $U \subset Y$ denote an open affine neighbourhood of $y$. Let
$X_U$ denote the inverse image of $U$ in $X$. By abuse of notation we
continue to denote by $f$ the induced map $X_U \to U$.
Let $i : \{ y \} \into U$ denote the inclusion of $\{ y \}$ and $j$
the inclusion of the complement $U \setminus \{ y \}$. In
the distinguished triangle
\[
j_!j^! f_* \RM_{X_U}[n] \to f_* \RM_{X_U}[n] \to i_*i^*f_* \RM_{X_U}[n] \triright
\]
all objects belong to ${}^pD^{\le 0}(U)$. Because $U$ is affine
$H^q(U, \FS) = 0$ if $q > 0$ for $\FS \in {}^pD^{\le 0}(U)$ by
Proposition \ref{prop:affinevanish}. In
particular
\[
r' : H^n(X_U) = H^0(U, f_* \RM_{X_U}[n]) \to H^0(U, i_*i^*f_*\RM_{X_U}[n])
= H^n(F)
\]
is surjective.

We may factor our map $r$ as $H^n(X) \to H^n(X_U) \stackrel{r'}{\to} H^n(F)$.
By mixed Hodge theory \cite[Prop. 8.2.6]{DH3} the images of $r$ and $r'$ agree. Hence $r$ is
surjective and the lemma follows.

\bigskip

We may now deduce the Semi-Small Index Theorem from
Theorem \ref{thm:lef}. We have an isometric embedding $cl : H_n(F) \into P^{m,m} \subset
H^n(X)$. By the Hodge-Riemann relations the Poincar\'e form on the
later space is $(-1)^m$-definite. Hence this is also the case for the intersection form
on $H_n(F)$.

\subsection{Hard Lefschetz via Positivity} \label{sec:hlp}

Our goal is to outline a proof of Theorem
\ref{thm:lef}, which we will carry out in the next section. 
Beforehand we recall an old idea
to prove the hard Lefschetz
theorem by combining Poincar\'e duality and the weak Lefschetz theorem with the Hodge-Riemann
relations in dimension one less.

To this end suppose that $X \subset \PM$ is a smooth projective
variety of dimension $n$ and let $D \subset X$ be a general (i.e. smooth) hyperplane
section. Consider the graded vector spaces
\begin{gather*}
  H = \bigoplus H^j \quad \text{where} \quad H^j := H^{n+j}(X),\\
  H_D = \bigoplus H_D^j \quad \text{where} \quad H_D^j := H^{n-1+j}(D).
\end{gather*}
In the following, we attempt to carry out an inductive proof of the
hard Lefschetz theorem for $H$. We assume as known the weak
Lefschetz theorem and Poincar\'e duality in general and the hard
Lefschetz theorem and Hodge-Riemann relations for $H_D$.

The inclusion $i : D \into X$ gives Poincar\'e dual restriction and Gysin morphisms
\[
i^* : H^j \to H_D^{j+1} \quad \text{and} \quad i_! : H_D^j \to H^{j+1}.
\]
Denote by $\omega$ the Chern class determined by our embedding $X \subset \PM$
and let $\omega_D$ denote its restriction to $D$. We have:
\begin{gather}
  \label{eq:12}
  \omega \wedge \alpha = i_! \circ i^* (\alpha) \quad \text{for all
    $\a \in H$,} \\
  \label{eq:13}
  \omega_D \wedge \b = i^* \circ i_! (\b) \quad \text{for all
    $\b \in H_D$.}
\end{gather}
Moreover, by the weak Lefschetz theorem:
\begin{gather}
  \label{eq:14}
  i^* : H^j \to H_D^{j+1} \text{ is an isomorphism if $j < -1$ and
    injective if $j = -1$,} \\
  \label{eq:15}
  i_! : H_D^{j-1} \to H^j \text{ is an isomorphism if $j > 1$ and
    surjective if $j = 1$.}
\end{gather}
Now the hard Lefschetz theorem for $H_D$ implies the hard Lefschetz
theorem for $\omega^k : H^{-k} \to H^k$ for $k >1$ because we can
factor $\omega^k$ as
\[
H^{-k} \simto H_D^{-k+1} \stackrel{\omega_D^{k-1}}{\longrightarrow} H_D^{k-1} \simto H^k
\]
where the first and last maps are weak Lefschetz isomorphisms.

The missing case is $\omega : H^{-1} \to H^1$. However in this case
one may use the relations \eqref{eq:12} and \eqref{eq:13} to deduce
that $i^*$ restricts to a map:
\[
i^* : P^{-1} = \ker ( \omega^2 : H^{-1} \to H^3) 
 \to P^0_D := \ker (\omega_D : H_D^0 \to H^2_D).
\]
Hence if $0 \ne \alpha \in P^{-1}$ is of pure Hodge type $(p,q)$ then, by
weak Lefschetz and the Hodge-Riemann bilinear relations,
\[
0 \ne (i^*\alpha, i^*\alpha) = (\alpha, \omega \wedge \alpha).
\]
It follows that $\omega : H^{-1} \to H^1$ is injective, and hence an isomorphism
(by Poincar\'e duality).

\begin{rema} \label{rem:deg0}
  The above line of reasoning can be used to deduce the
  Hodge-Riemann bilinear relations for all primitive subspaces $P^j \subset H^j$ with $j <
  0$. However the crucial case of the Hodge-Riemann relations for the
  middle degree $P^0 \subset H^0$ is missing. Hence we cannot close the induction.
\end{rema}

\subsection{Hard Lefschetz and Hodge-Riemann for Semi-Small Classes}

We now outline de Cataldo and Migliorini's proof of Theorem
\ref{thm:lef}. The basic idea is to combine the argument of the
previous section with a limit argument to recover the missing
Hodge-Riemann relations. Recall that
\[
f : X \to Y
\]
is a semi-small morphism with
$X$ connected, smooth and projective. The proof is by induction on
the dimension $n$ of $X$. If $n = 0, 1$ then $f$ is finite, and the
theorem can be checked by hand.

\emph{Step 1: Hard Lefschetz.} Let $\LC$ be an ample line bundle on
$Y$, $i : D \into Y$ the inclusion of a general hyperplane section, $f_D : X_D :=
f^{-1}(D) \to D$ the induced map, $\omega \in H^2(X)$ the Chern class
of $f^* \LC$, and $\omega_D$ its restriction of $X_D$.

A Bertini type argument (see \cite[Prop. 2.1.7]{dCM}) guarantees that
$X_D$ is smooth and that $f_D$ is semi-small. Hence we can apply
induction to deduce that hard Lefschetz and the Hodge-Riemann relations hold for the action of $\omega_D$ on
$H^*(X_D)$. Because $f_* \RM_X[n]$ is perverse, the weak Lefschetz
theorem holds for the restriction map
\[
i^* : H^{*+n}(X) = H^*(Y, f_* \RM_X[n]) \to H^*(Y, i_*i^*f_*\RM_X[n]) = H^{*+n}(X_D)
\]
and its dual. Now the arguments of the previous section allow us to deduce
that $\omega$ satisfies hard Lefschetz on $H^*(X)$.

\emph{Step 2: Hodge-Riemann.} We explain how to deduce the
Hodge-Riemann relations for the crucial case $H^0 =
H^n(X)$. Hodge-Riemann relations in degrees $< 0$ follow similarly
(or alternatively one can use Step 1 and Remark \ref{rem:deg0}).

Let $\eta$ denote an ample class on $X$. Then $\omega + \e \eta$
belongs to the ample cone for all $\e > 0$. For $\e \ge 0$ consider the subspaces:
\begin{gather*}
  P_\e^0:= \ker( (\omega + \e\eta) : H^{0} \to H^{2} ), \\
  P_\e^{p,q} := (P_\e^0)_{\CM} \cap H^{p,q} \quad \text{where $p +
    q  = n$.}
\end{gather*}

We claim that, in the Grassmannian of subspaces of $H^0$, we have
\begin{equation}
\lim_{\e \to 0} P_\e^{p,q} = P_0^{p,q}. \label{eq:lim}
\end{equation}
The left hand side is clearly contained in the right hand side. The
claim now follows because both sides have dimension $\dim H^{p,q} -
\dim H^{p+1,q+1}$ (for the left hand side this follows via classical Hodge
theory and Remark \ref{rem:amplecone}, for the right hand side it follows by hard Lefschetz for
$\omega$ established in Step 1).

Recall our Hermitian form $\kappa(\a, \b) = (\sqrt{-1})^n \int \a
\wedge \overline{\b}$ on $H^0_{\CM}$. 
We conclude from \eqref{eq:lim} that any $\a \in P_0^{p,q}$ is a limit of classes in
$P^{p,q}_\e$ as $\e \to 0$. Hence, by the Hodge-Riemann relations for
the classes $\omega + \e\eta$ (which lie in the ample cone) we have
\begin{equation}
(\sqrt{-1})^{p-q-n}(-1)^{n(n-1)/2} \kappa(\a,\overline{\a}) \ge
0 \quad \text{for any $\a \in P^{p,q}_0$.}\label{eq:15}
\end{equation}
By Hard Lefschetz the restriction of $\kappa$ to each $P^{p,q}_0$ is
non-degenerate. However \eqref{eq:15} tells us that our Hermitian form is also
semi-definite on $P^{p,q}_{0}$. We conclude that our form is
definite and we have a strict inequality
\begin{equation}
(\sqrt{-1})^{p-q-n}(-1)^{n(n-1)/2} \kappa(\a,\overline{\a}) >
0 \quad \text{for any $\a \in P^{p,q}_0$.}
\end{equation}
This yields the Hodge-Riemann relations for $H^0$.

\section{General maps}\label{sec:general}

In this section we outline de Cataldo and Migliorini's proof of the
Decomposition Theorem for general projective maps. The proof follows
the same main lines as the semi-small case, however the collection of
statements needed through the 
induction is more involved. We
refer to this collection as the \og Decomposition
Theorem Package\fg{}. We begin by stating all theorems constituting
the package, and then proceed to an outline of the inductive proof.

\subsection{The Decomposition Theorem Package}

We assume as always that $X$ is a smooth connected projective variety of
dimension $n$ and that
\[
f : X \to Y
\]
is a surjective projective morphism. Let us fix the following
two classes in $H^2(X)$:
\begin{gather*}
  \eta : \text{the Chern class of a relatively ample (with respect to $f$) line bundle,} \\
  \beta : \text{the Chern class of the pull-back (via $f$)
    of an ample line bundle on $Y$.}
\end{gather*}

Because $H^2(X) = \Hom_{D^b_c(X)}(\RM_X, \RM_X[2])$ we may interpret
$\eta$ as a map $\eta : \RM_X[n] \to \RM_X[n+2]$. Pushing forward we
obtain a map (also denoted $\eta$):
\[
\eta : f_* \RM_X[n] \to f_*\RM_X[n+2].
\]
Recall that every object in $D^b_c(Y)$ carries a perverse
filtration. Moreover this filtration is preserved by any morphism in
$D^b_c(Y)$. Thus $\eta$ induces maps (for all $m \in \ZM$):
\begin{equation}
\eta : \ptau_{\le m} f_*\RM_X[n] \to \ptau_{\le m+2}
f_*\RM_X[n].\label{eq:ptaueta}
\end{equation}

The Relative Hard Lefschetz Theorem concerns the associated graded of $\eta$:

\begin{theo}[Relative Hard Lefschetz Theorem] \label{thm:rh}
For $i \ge 0$, $\eta$ induces an isomorphism
\[
\eta^i : \pH^{-i}(f_*\RM_X[n]) \simto \pH^i(f_*\RM_X[n]).
\]
\end{theo}

\begin{rema}
 The Relative Hard Lefschetz Theorem specialises to the Hard Lefschetz
 Theorem if $Y$ is a point. If $f$ is smooth (the setting of Deligne's
 Theorem) the Relative Hard Lefschetz Theorem follows from the
 classical Hard Lefschetz Theorem applied to the fibres of $f$. If $f$
 is semi-small then $\pH^{-i}(f_*\RM_X[n]) = 0$ unless $i = 0$ and the
 Relative Hard Lefschetz Theorem holds trivially.
\end{rema}

It is a formal consequence of the Relative Hard Lefschetz Theorem that
we have a decomposition (see \cite{DD}) 
\begin{equation} \label{eq:rhlsplitting}
f_*\RM_X[n] \cong \bigoplus \pH^i(f_*\RM_X[n])[-i].
\end{equation}
The heart of the Decomposition Theorem is now:

\begin{theo}[Semi-Simplicity Theorem] \label{thm:ss}
Each $\pH^i(f_*\RM_X[n])$ is a semi-simple perverse sheaf.
\end{theo}

\begin{rema} If $f$ is smooth (the setting of Deligne's Theorem) the
  Semi-Simplicity Theorem follows from the fact (\cite[Theorem 7.1]{G}, \cite[Theorem
  4.2.6]{DH2}, \cite[Theorem 7.25]{S})
  that a local system underlying a polarisable pure variation of Hodge
  structure on a smooth variety is
  semi-simple. If $f$ is semi-small then all the content of the
  Decomposition Theorem is contained in the Semi-Simplicity Theorem
  for $\pH^{0}(f_*\RM_X[n]) = f_* \RM_X[n]$.
\end{rema}

\begin{rema} \label{rem:fakefibres}
By the Semi-Simplicity Theorem, we have a
canonical isomorphism
\[
\pH^i(f_*\RM_X[n]) = \bigoplus V_{\l, \LS,i} \otimes
IC(\overline{Y}_\l, \LS)
\]
where the sum runs over all pairs $(Y_\l, \LS)$ consisting of a
stratum $Y_\l$ and an (isomorphism class of)
simple local system $\LS$ on $Y_\l$, and $V_{\l,\LS,i}$ is a real
vector space. By semi-simplicity the map 
\[
\eta : \pH^i(f_*\RM_X[n]) \to \pH^{i+2}(f_*\RM_X[n])
\] is completely described by maps of vector spaces $\eta : V_{\l,\LS,i} \to
V_{\l,\LS,i+2}$ for all pairs $(Y_\l,\LS)$. The Relative Hard Lefschetz theorem now becomes the
statement that the degree two endomorphism $\eta$ of the finite-dimensional graded vector space
\[
V_{\l, \LS} := \bigoplus V_{\l, \LS,i}
\]
satisfies hard Lefschetz for all pairs $(Y_\l,\LS)$. 
\end{rema}

As in the semi-small case it is important to understand the
structure that the above theorems give on the global cohomology of $X$. We set
$H^i := H^{n+i}(X) = H^i(Y, f_*\RM_X[n])$ as usual. By taking global cohomology of the 
perverse filtration
\[
\dots \to \ptau_{\le m} f_* \RM_X[n] \to \ptau_{\le m+1} f_* \RM_X[n]
\to \dots
\]
we obtain the (global) \emph{perverse filtration} on $H$:
\[
\dots \subset H_{\le m} \subset H_{\le m+1} \subset \dots
\]
Recall that $H$ is equipped with its Poincar\'e
form. With respect to this form one has
\begin{equation}
H_{\le i}^{\perp} = H_{< -i}.\label{eq:19}
\end{equation}

Consider the associated graded of the perverse filtration:
\[
H_i := H_{\le i}/H_{<i} \quad \text{and} \quad \gr H = \bigoplus H_i.
\]
By \eqref{eq:19} the Poincar\'e form induces a non-degenerate form
$H_i \times H_{-i} \to \RM$ and hence a non-degenerate form on $\gr H$.

\begin{prop} \label{prop:hs}
  The perverse filtration is a filtration by pure Hodge
  substructures. In particular, each $H_i^j$ is a pure Hodge structure of
  weight $n+j$.
\end{prop}

The action of $\eta$ and $\beta$ on $H$ satisfy:
\begin{gather}
  \label{eq:beta}
  \b(H_{\le m}) \subset H_{\le m} \quad \text{for all $m \in \ZM$,} \\
  \label{eq:eta}
  \eta(H_{\le m}) \subset H_{\le m+2} \quad \text{for all $m \in \ZM.$}
\end{gather}
(The first inclusion follows
because the cohomology of any complex on $Y$ is a graded module over
$H^*(Y)$. The second inclusion follows from \eqref{eq:ptaueta}.) Hence
we obtain operators
\[
\b : H_i^j \to H_i^{j+2} \quad \text{and} \quad \eta : H_i^j
\to H_{i+2}^{j+2}.
\]

\begin{theo}[Relative Hard Lefschetz in Cohomology] \label{thm:rhL}
For $i \ge 0$, $\eta$ induces an isomorphism
\begin{equation} \label{eq:relcoh}
\eta^i : H_{-i} \simto H_{i}.
\end{equation}
\end{theo}

\begin{theo}[Hard Lefschetz for Perverse Cohomology] \label{thm:phL}
For all $i \in \ZM$ and $j \ge 0$, $\b$ induces an isomorphism
\begin{equation} \label{eq:hlp}
\beta^j : H_i^{i-j}\simto H_{i}^{i+j}.
\end{equation}
\end{theo}


\begin{rema}
  One may depict the $H_i^j$ and maps $\beta$ and
$\eta$ as a two-dimensional array:
\begin{equation} \label{eq:array}
\begin{array}{c}
  \begin{tikzpicture}[scale=1.4]
 \node (z-2-2) at (-2,0) {$H_{-2}^{-2}$};
\node (z-1-2) at (-1,-1) {$H_{-1}^{-2}$};
\node (z-1-1) at (-1,0) {$H_{-1}^{-1}$};
\node (z-10) at (-1,1) {$H_{-1}^{0}$};
\node (z02) at (0,2) {$H_0^2$};
\node (z01) at (0,1) {$H_0^1$};
\node (z00) at (0,0) {$H_0^0$};
\node (z0-1) at (0,-1) {$H_0^{-1}$};
\node (z0-2) at (0,-2) {$H_0^{-2}$};
\node (z12) at (1,1) {$H_{1}^{2}$};
\node (z11) at (1,0) {$H_{1}^{1}$};
\node (z10) at (1,-1) {$H_{1}^{0}$};
 \node (z22) at (2,0) {$H_{2}^{2}$};
%
\node (z-11) at (-1.5,1.5) {$\ddots$};
\node (z1-1) at (1.5,-1.5) {$\ddots$};
\node[xscale=-1,yscale=1] at (1.5,1.5) {$\ddots$}; 
\node[xscale=-1,yscale=1] at (-1.5,-1.5) {$\ddots$}; 
\draw[->] (z00) .. controls (1,0.4) .. node[above, near end] {\tiny $\eta$} (z22);
\draw[->] (z-2-2) .. controls (-1,0.4) .. node[above, near end] {\tiny
  $\eta$} (z00);
\draw[->] (z00) ..controls (-0.5,1) .. node[left, near end] {\tiny $\beta$} (z02);
\draw[->] (z0-2) ..controls (-0.5,-1) .. node[left, near end] {\tiny $\beta$}
(z00);
\draw[->] (-2.9,-1) to node[above,sloped] {\small $j-i$} (-2.9,1);
\draw[->] (-1,-2.6) to node[below,sloped] {\small $i$} (1,-2.6);
%
%
\end{tikzpicture}
\end{array}
\end{equation}
(We have only depicted the maps with source or target $H_0^0$.)
Relative Hard Lefschetz states that $\eta$ satisfies Hard Lefschetz
along each row, and Hard Lefschetz for Perverse Cohomology states that
$\b$ satisfies Hard Lefschetz along each column.
\end{rema}

Recall that the hard Lefschetz theorem leads to a primitive
decomposition of cohomology. The above two theorems lead to
a bigraded primitive decomposition; set
\[
P_{-i}^{-j} := \ker (\eta^{i+1} : H_{-i}^{-i-j} \to H_{i+2}^{i-j +
  2}) \cap \ker (\b^{j+1} : H_{-i}^{-i-j} \to H_{-i}^{-i + j +
  2}) \subset H_{-i}^{-i-j}.
\]

\begin{coro}[$(\eta,\beta)$-Primitive Decomposition] \label{cor:eb}
The inclusions
  $P_{-i}^{-j} \into H_{-i}^{-i-j}$ induce a 
  canonical isomorphism of $\RM[\eta, \beta]$-modules:
\[
\bigoplus_{i, j \ge 0} \RM[\eta]/(\eta^{i+1})\otimes \RM[\beta]/(\beta^{j+1})
\otimes P_{-i}^{-i-j} \simto \gr H.
\]
\end{coro}

\begin{rema} Recall that the Hard Lefschetz Theorem can be rephrased
  in terms of an $\mathfrak{sl}_2$-action (see Remark
  \ref{rem:sl2}). Similarly, Theorems \ref{thm:rhL} and \ref{thm:phL} are equivalent to the
  existence of an $\mathfrak{sl}_2 \times \mathfrak{sl}_2$-action on
  $\gr H$ such that, for all $x \in H_i^j$, we have
\[
e_1(x) = \eta(x), \quad h_1(x) = i(x), \quad e_2(x) = \beta(x), \quad h_2(x) = (j-i)(x)
\]
(the subscript indicates in which copy of $\mathfrak{sl}_2$ the
generator lives). The $(\eta, \beta)$-decomposition coincides with the isotypic
decomposition and the primitive subspaces $P_{-i}^{-j}$ are the lowest weight spaces.
\end{rema}

Our form on $\gr H$ induces a form on each $H_{-i}^{-i-j}$ for
$i , j \ge 0$; if $[\a], [\b]
\in H_{-i}^{-i-j}$ are represented by classes $\a, \b \in H_{\le -i}^{-i-j}$ we
set
\[
S_{ij}(\a, \b) := \int_X \eta^i \wedge \beta^j \wedge \a \wedge \b.
\]
This form is well defined and non-degenerate.

\begin{theo}[$(\eta,\beta)$-Hodge-Riemann Bilinear Relations] \label{thm:ebhr}
Each Hodge structure $P_{-i}^{-j}$ (of weight $n - i - j$) is polarised by the forms $S_{ij}$.
\end{theo}

\begin{rema}
  With the notation in Remark \ref{rem:fakefibres} we have
  \begin{equation}
\gr H = \bigoplus V_{\l, \LS} \otimes
IH(\overline{Y}_\lambda,\LS).\label{eq:grdecomp}
\end{equation}
In the array \eqref{eq:array} the columns (resp. rows) correspond to
the grading on 
$IH(\overline{Y}_\lambda,\LS)$
 (resp.  $V_{\l,
  \LS}$). The above theorems can be understood as saying that \og each
row and column looks like the cohomology of a
smooth projective variety\fg{}. This remarkable point of view is
emphasised in \cite{Mac}. See \cite{GM} and \cite[\S 2]{dCM2} for examples.
\end{rema}

\begin{rema} Any choice of an isomorphism of complexes as in
  \eqref{eq:rhlsplitting} gives
  an isomorphism $H \cong \gr H$ of vector spaces. It is possible to
  choose this isomorphism so as to obtain an isomorphism of Hodge
  structures \cite{dCM7}. Under
  such an isomorphism the decomposition $H = \bigoplus V_{\l, \LS} \otimes
IH(\overline{Y}_\lambda,\LS)$ given by \eqref{eq:grdecomp} is
of motivic nature \cite{dCM5}; that is, the projectors in this decomposition are motivated cycles in the sense of
Andr\'e \cite{A}, and are given by algebraic cycles if
Grothendieck's standard conjectures are true.
\end{rema}

\subsection{Defect of Semi-Smallness and Structure of the Proof}

We now give an outline of the argument. For a projective
and surjective map $f : X \to Y$ with $\dim X = n$ its \emph{defect of
  semi-smallness} is
\[
r(f) := \textrm{max} \{ i \in \ZM \; | \; \pH^i(f_*\RM_X[n]) \ne 0\}.
\]
Equivalently, if $Y := \bigsqcup Y_\lambda$ is a stratification of
$f$ and we choose a point $y_\l \in Y_\l$ in each stratum then
\[
r(f) = \textrm{max}_{\l \in \Lambda} \{  \; 2 \dim f^{-1}(y_\lambda) +
\dim Y_\lambda - \dim X \}.
\]
We have $r(f) \ge 0$ and $r(f) = 0$ if and only if $f$ is semi-small.

The proof is via simultaneous induction on the defect of
semi-smallness and on the dimension of the image of $f$. That is, if we
fix $f$ we may assume that the Decomposition Theorem Package is known for any
projective map $g : X' \to Y'$ with $r(g) < r(f)$ or $r(g) = r(f)$ and
$\dim g(X') < \dim f(X)$. The base case is when $f$ is the projection
to a point, in which case all statements follow from classical Hodge
theory.

The proof breaks up into four main steps. The titles in italics are
those of the upcoming sections. Only steps 2) and 3) have an analogue in the case
of a semi-small map:

\begin{enumerate}
\item \emph{Relative Hard Lefschetz via Semi-Simplicity.} The Relative
  Hard Lefschetz Theorem is 
  deduced from the Decomposition Theorem for the relative universal hyperplane
  section morphism. The decomposition
\[
f_*\RM_X[n] \cong \bigoplus \pH^i(f_*\RM_X[n])[-i].
\]
is an immediate consequence, as is the semi-simplicity of
$\pH^i(f_*\RM_X[n])$ for $i \ne 0$.
\item \emph{Miraculous Approximability.} All global statements in the Decomposition
  Theorem Package are established. The crucial case $P_0^0$ is
  established by a limit argument.
\item \emph{Local Study of the Decomposition Theorem.} Analogously to the case of a semi-small map,
  the $(\eta,\beta)$-Hodge-Riemann Bilinear Relations are used to show that we have a splitting
\begin{equation} \label{eq:split}
\pH^0(f_*\RM_X[n]) = \bigoplus IC(\overline{Y}_\l, \LS_\l)
\end{equation}
for certain local systems $\LS_\l$.
\item \emph{Semi-Simplicity of Local Systems.} Deligne's theorem is
  used to show that each local system
  $\LS_\lambda$ appearing in \eqref{eq:split} is semi-simple.
\end{enumerate}

\begin{rema}
After the first step above the Semi-Simplicity
Theorem is easily reduced to proving the semi-simplicity of the middle primitive summand
\[
\mathcal{P}^0 := \ker (\eta : \pH^0(f_*\RM[n]) \to \pH^2(f_*\RM[n])).
\]
Philosophically, the proof of the general case should  involve
simply repeating the proof of the semi-small case for the summand
$\mathcal{P}^0$. This point of view is explained  in \cite[\S 3.3.2]{dCM4} where de
Cataldo and Migliorini call $\mathcal{P}^0$ the \og semi-small
soul\fg{} of the map $f$.
\end{rema}

\begin{rema} In the approach of Beilinson-Bernstein-Deligne-Gabber and
  Saito the Relative Hard Lefschetz Theorem and Decomposition Theorem
  are deduced from purity; in their approach one
  gets the Decomposition Theorem \og all at once\fg{}. The situation is quite different in
  de Cataldo and Migliorini's proof, where the Relative Hard Lefschetz
  Theorem (deduced from the Semi-Simplicity Theorem for a map with
  smaller defect of semi-smallness) is a crucial stepping stone in the induction.
\end{rema}

\subsection{Relative Hard Lefschetz via Semi-Simplicity}

We have explained in \S \ref{sec:hlp} how the Hard Lefschetz Theorem and the Hodge-Riemann 
Relations in dimensions $\le n-1$ imply the hard Lefschetz theorem in
dimension $n$. This step relies on positivity in a crucial way. In
this section we explain an older approach (due to Lefschetz) which
deduces hard Lefschetz from a semi-simplicity
statement. This approach is used to
prove the Relative Hard Lefschetz Theorem.

We first recall the idea in the absolute case. Suppose that $X \subset
\PM$ is a smooth projective variety of dimension $n$. We set $H^i :=
H^{n+i}(X)$ and let $\eta \in H^2(X)$ denote the Chern class
of the embedding. We
have explained in \S~\ref{sec:hlp} why weak Lefschetz and the Hard
Lefschetz Theorem for a smooth hyperplane section allow us to deduce
hard Lefschetz for $H$, except for the crucial ``middle'' case:
\begin{equation}
  \label{eq:crucial}
  \eta : H^{-1} \to H^{1}.
\end{equation}
If $D \subset X$ denotes a smooth hyperplane section, we also
explained that we may factor \eqref{eq:crucial} as the composition of the
restriction and its dual:
\begin{equation}
\eta  : H^{-1} \stackrel{i^*}{\into} H^{n-1}(D)
\stackrel{i_!}{\onto} H^{1}.\label{eq:factor}
\end{equation}
We now give a geometric description of the image of $i^*$.

Let $\PM^\vee$ denote the complete linear system of hyperplane
sections of $X$ and let $\YC \subset \PM^\vee$ denote the open subvariety of
smooth hyperplane sections. The morphism
\[
g : \XC := \{ (x, s) \in X \times \YC \;|\; s(x) = 0 \} \to \YC
\]
induced by the projection is the \emph{universal hyperplane
  section morphism}. Its fibres are the smooth hyperplane sections of
$X$. The map $g$ is smooth and proper and
\[
\LS := R^{n-1}g_* \RM_X
\]
is a local system whose fibre at $D' \in \YC$ is
$H^{n-1}(D')$.

Recall our chosen hyperplane $D$ from above. Regarding $D \in
\YC$ as a basepoint, we may alternately view $\LS$ as providing us
with a representation of the fundmental group $\pi_1 := \pi_1(\YC,D)$ on
$H^{n-1}(D)$. The following two results are fundamental observations
of Lefschetz.  For a modern proof see \cite[\S 4]{DW2}.

\begin{prop}  \label{prop:LefschetzImage}
We have a commutative diagram:
  \begin{equation*}
    \begin{tikzpicture}[xscale=1.3]
      \node (h-1) at (-2,0) {$H^{-1}$};
      \node (h0) at (0,0) {$H^{n-1}(D)$};
      \node (h1) at (2,0) {$H^{1}$};
      \node (inv) at (-1,-1) {$H^{n-1}(D)^{\pi_1}$};
      \node (coinv) at (1,-1) {$H^{n-1}(D)_{\pi_1}$};
\draw[right hook-latex] (h-1) to node[above] {$i^*$} (h0);
\draw[->] (h-1) to node[sloped,above] {$\sim$} (inv);
\draw[->] (coinv) to node[sloped,above] {$\sim$} (h1);
\draw[->>] (h0) to node[above] {$i_!$} (h1);
\draw[right hook-latex] (inv) to (h0);
\draw[->>] (h0) to (coinv);
    \end{tikzpicture}
  \end{equation*}
(where $V^{\pi_1}$ and $V_{\pi_1}$ denotes $\pi_1$-invariants and
coinvariants respectively).
\end{prop}

 \begin{coro}
   If $H^{n-1}(D)$ is semi-simple as a $\pi_1$-module then $\eta : H^{-1}
   \to H^1$ is an isomorphism.
 \end{coro}


We now return to our setting of $f : X \to Y$ a projective morphism
with $X$ smooth and projective of dimension $n$. For simplicity we fix
an embedding $X \subset \PM$ of $X$ into a projective space of
dimension $d$ and let $\PM^\vee$ denote the dual projective
space. Consider the following spaces and maps:
\begin{gather*}
  \begin{tikzpicture}[scale=2]
    \node (X) at (-1,1) {$X$};
 \node (X1) at (1,1) {$X \times \PM^\vee$};
 \node (Xt) at (4,1) {$\XC := \{ (x,s) \in X \times \PM^\vee \; | \;
   s(x) = 0 \}$};
 \node (p) at (-1,0) {$Y$};
 \node (P1) at (1,0) {$\YC := Y \times \PM^\vee$};
\draw[->] (X1) to node[above] {$p$} (X);
\draw[->] (P1) to node[above] {$p$} (p);
\draw[->] (X) to node[left] {$f$} (p);
\draw[->] (X1) to node[right] {$f$} (P1);
\draw[->] (Xt) to node[above] {$i$} (X1);
\draw[->] (Xt) to node[below] {$g$} (P1);
  \end{tikzpicture}
\end{gather*}
(Note that different arrows have the same name.) The map $g$ is
the \emph{(relative) universal hyperplane section morphism}. Its fibre
over a point $(y, s) \in \YC$ is the intersection of $f^{-1}(y)$ with the
hyperplane section of $X$ determined by $s$.

The following crucial lemma (\og the defect of semi-smallness goes
down\fg{}) allows us to apply our inductive
assumptions to conclude that the Decomposition Theorem Package holds
for $g$:

\begin{lemm}If $r(f) > 0$ then $r(g) < r(f)$. If $r(f) = 0$ then $r(g)
  = 0$.
\end{lemm}

The proof is an easy analysis of a stratification of $g$, see
\cite[Lemma 4.7.4]{dCM}.

Set $m := n + d - 1 = \dim \XC$.

\begin{prop}[Relative Weak Lefschetz for Perverse Sheaves] \label{prop:rwl}
  \begin{enumerate}
  \item  For $j < -1$ there is a natural isomorphism: 
\[
p^*( \pH^j(f_* \RM_X[n]))[d] = \pH^j(f_*\RM_{X \times \PM^\vee}[m+1]) \simto \pH^{j+1}(g_*\RM_{\XC}[m]).
\]
  \item  For $j > 1$ there is a natural isomorphism: 
\[
\pH^{j-1}(g_*\RM_{\XC}[m]) \simto \pH^j(f_*\RM_{X \times \PM^\vee}[m+1]) = p^*( \pH^j(f_* \RM_X[n]))[d].
\]
\item $p^*(\pH^{-1}(f_*\RM_X[n]))[d]$ is the largest 
  subobject of $\pH^0(g_*\RM_{\XC}[m])$ coming from $Y$.
\item $p^*(\pH^{1}(f_*\RM_X[n]))[d]$ is the largest quotient of
  $\pH^0(g_*\RM_{\XC}[m])$ coming from $Y$.
  \end{enumerate}
\end{prop}

For a proof of these statements, see \cite[5.4.11]{BBD}. The first two
statements are a consequence of the fact that the restriction of
$f$ to the complement of $\XC \subset X \times \PM^\vee$ is affine,
combined with the cohomological dimension of affine morphisms. The
second two statements are relative analogues of Proposition
\ref{prop:LefschetzImage}. (The notion of largest subobject or
quotient coming from $Y$ is well defined because $p$ is a smooth
morphism with connected fibres. Thus $p^*[d]$ identifies the category of
perverse sheaves on $Y$ with a full subcategory of perverse
sheaves on $\YC$, see \cite[\S 4.2.5--6]{BBD}.)

\begin{rema}
  The semi-simplicity of $\pH^j(f_*\RM_X)$ for $j \ne -1, 0, 1$ is an
  immediate consequence of (1) and (2) above. The
  semi-simplicity of $\pH^j(f_*\RM_X)$ for $j = -1, 1$ follows from
  (3) and (4).
\end{rema}

Let us now explain the proof of the Relative Hard Lefschetz Theorem,
following \cite[\S 5.4.10]{BBD}. The adjunction morphism 
\[
\RM_{X\times \PM^\vee}[m+1] \to i_*i^*\RM_{X\times \PM^\vee}[m+1] = i_*\RM_\XC[m+1]
\]
induces morphisms
\[
i^* : \pH^j(f_*\RM_{X\times \PM^\vee}[m+1]) \to \pH^{j+1}(f_*\RM_\XC[m]).
\]
For $j < 0$ these are the morphisms appearing in parts (1) and (3) of
Proposition \ref{prop:rwl}. Taking duals we obtain morphisms (we use that $\XC
\subset X \times \PM^\vee$ is smooth)
\[
i_! : \pH^{j-1}(f_*\RM_\XC[m]) \to \pH^j(f_*\RM_{X\times \PM^\vee}[m+1]).
\]
For $j > 0$ these are the morphisms appearing in (2) and (4) of
Proposition \ref{prop:rwl}.

For $j \ge 0$ consider the morphisms:
\[
\pH^{-j-1}(f_*\RM_{X\times \PM^\vee}[m+1]) \stackrel{i^*}{\to} \pH^{-j}(f_*\RM_\XC[m])
\stackrel{\eta^j}{\to} \pH^j(f_*\RM_\XC[m]) \stackrel{i_!}{\to} \pH^{j+1}(f_*\RM_{X\times \PM^\vee}[m+1]) 
\]
We claim that for $j \ge 0$ the composition is an isomorphism. If $j >
0$ this follows because the first and last maps are isomorphisms by
Proposition \ref{prop:rwl} and the middle map is an isomorphism by
the Relative Hard Lefschetz Theorem for $g$. For $j = 0$ the
composition is an isomorphism by parts (3) and (4) of
Proposition \ref{prop:rwl} and the Semi-Simplicity Theorem for
$\pH^{0}(f_*\RM_\XC[m])$ (which is known by induction).

Finally, the above composition agrees up to shift with the pull-back
via $p$ of the map
\[
\eta^j : \pH^{-j-1}( f_*\RM_X[n]) \to \pH^{j+1}( f_*\RM_X[n]).
\]
Hence $\eta^j$ is an isomorphism for $j \ge 0$ and the Relative Hard
Lefschetz Theorem follows.

\begin{rema} \label{rem:hlpn}
One may also use Proposition \ref{prop:rwl} and induction to deduce the Hard
Lefschetz Theorem for Perverse Cohomology (Theorem \ref{thm:phL}) for all
degrees except perverse degree zero. Indeed, if $\FS$ is a complex of
sheaves on $Y$ and $\beta \in H^2(Y)$ is an ample class, then $\beta$
satisfies hard Lefschetz on $H(Y, \FS)$ if and only if $\b  + \z \in
H^2(Y \times \PM^\vee)$ satisfies hard Lefschetz on $H(Y \times
\PM^\vee, p^*\FS[d]) = H(Y, \FS) \otimes H(\PM^\vee)$, where $\z$ denotes the pull-back of any
non-zero element of $H^2(\PM^\vee)$.
\end{rema}

\subsection{Miraculous Approximability}

In this section we discuss the inductive proofs of the global
statements in the Decomposition Theorem package. Here the proofs are
often routine and sometimes technical and we will not attempt to give more than a rough
outline. For more detail the reader is referred to \cite[\S
5.2-5.4]{dCM2}.

What do we know at this stage? The Relative Hard Lefschetz Theorem
proved in the previous step implies immediately the Relative Hard
Lefschetz Theorem in Cohomology (Theorem \ref{thm:rhL}). Also, the
previous step gives the Hard Lefschetz Theorem in Perverse Cohomology
(Theorem \ref{thm:phL}) except for perverse degree zero (i.e. 
$i = 0$ in Theorem \ref{thm:phL}) by Remark \ref{rem:hlpn}.

As in the semi-small case, an argument involving a generic hyperplane section
$D \subset Y$ and the $(\eta,\beta)$-Hodge-Riemann relations for the
restriction of $f$ to $f^{-1}(D)$ 
allows us to deduce the missing $i = 0$ case of Theorem
\ref{thm:phL}. (This technique could also be used to prove
Theorem \ref{thm:phL} in the other cases, and avoid Remark
\ref{rem:hlpn}.) Theorem
\ref{thm:rhL} and Theorem \ref{thm:phL} and 
some linear algebra imply the $(\eta, \beta)$-Primitive Decomposition
(Corollary \ref{cor:eb}).

All that remains are the $(\eta,\beta)$-Hodge-Riemann Bilinear Relations.
To make sense of these relations we need Proposition \ref{prop:hs}, which tells us that the perverse
filtration on $H$ and its subquotients are pure Hodge structures. In
de Cataldo and Migliorini's original proof this fact was deduced from
Theorem \ref{thm:phL}. The idea is that one can use hard
Lefschetz on each $H_i$ to give a purely linear algebraic
definition of the perverse filtration (as a \og weight filtration\fg{}
associated to the operator $\b$), which then implies that it is
linear algebraic in nature, and hence is a filtration by pure Hodge structures.

However a more recent theorem of de Cataldo and Migliorini \cite{dCM3}
gives a conceptually and practically superior proof of Proposition
\ref{prop:hs}. Their result is that the perverse filtration of any
complex on $Y$ is given (up to shift) by a \og flag filtration\fg{}
associated to general flag of closed subvarieties of
$Y$. 
As a consequence the perverse filtration associated to any map is
by mixed Hodge structures (this result is independent of the
Decomposition Theorem and even holds over the
integers). Proposition \ref{prop:hs} is an easy consequence.

It remains to discuss the proof of the $(\eta,
\beta)$-Hodge-Riemann relations (Theorem \ref{thm:ebhr}). By taking
hyperplane sections in $X$ one may deduce the
$(\eta,\beta)$-Hodge-Riemann Bilinear Relations 
for the primitive subspaces $P_{-i}^{-j} \subset H_{-i}^{-i-j}$ for
all $i, j \ge 0$ with $(i, j) \ne (0,0)$. This reduction is
formally analogous to the semi-small case.

In the semi-small case the missing Hodge-Riemann relations in
degree 0 were deduced via a limit argument. Here the approach is
similar but more involved. The complication is that $P_0^0$ is
a subquotient of $H$, and so it is a priori not clear how to
realise $P_0^0$ as a \og limit\fg{} of a subspace in $H$. That this is
still possible explains the title of this section.

We proceed as follows. For $\e > 0$, $\b + \e \eta \in H^2(X)$
lies in the ample cone. Hence if
\[
\Lambda_\e := \ker( \b + \e \eta : H^0 \to H^2) \subset H^0 = H^n(X)
\]
then $d :=\dim \Lambda_\e = \dim H^0 - \dim H^2$ by hard Lefschetz for
$\b + \e \eta$ (see Remark \ref{rem:amplecone}). Consider the limit
(taken in the Grassmannian of $d$-dimensional
subspaces of $H^0$):
\[
\Lambda := \lim_{\e \to 0} \Lambda_\e.
\]

Each $\Lambda_\e$ is a polarised Hodge
substructure of $H^0$. Hence $\Lambda$ is a Hodge substructure
(being a Hodge substructure is a closed condition). Also, $\Lambda$ is
semi-polarised (i.e. the Hodge-Riemann relations hold for
$\Lambda$ if we replace strict inequality $> 0$ by $\ge 0$).

To keep track of degrees let us denote the map $\b : H^{-i} \to
H^{-i+2}$ by $\b_i$. Of course $\Lambda \subset \ker( \b_0)$ however 
equality does not hold in general because
\begin{equation} \label{eq:dimb0}
\dim \ker(\b_0) = \dim \Lambda + \dim \ker(b_2)
\end{equation}
as follows, for example, by noticing that one can perform this
calculation on $\gr H$.

It is important to be able to identify $\Lambda \subset H^0$
intrinsically. This is completed 
in \cite[\S 5.4]{dCM2}. As a consequence they deduce:

\begin{lemm}\label{lem:split} We have an orthogonal decomposition $\ker \b_0 = \Lambda \oplus \eta( \ker \b_2).$
\end{lemm}

The Poincar\'e form on the image of $\ker \b_0$ in $H_0^0$ is
non-degenerate by the Hard Lefschetz Theorem for $H_0$. Thus the
radical of the Poincar\'e form on $\ker \b_0$ is 
$\ker \b_0 \cap H^0_{<0}$. It follows from Lemma \ref{lem:split} that
the radical of the Poincar\'e form on $\Lambda$ is
\[
\Lambda_{<0} := \Lambda \cap H^0_{<0}.
\]
Thus the Poincar\'e form on
\[
\Lambda_0 := \Lambda / \Lambda_{<0}
\]
is a non-degenerate semi-polarisation, and thus a polarisation.

Finally, with a little more work the above lemma also implies
that we have an embedding of Hodge structures
\[
P_0^0 \into \Lambda_0
\]
which proves the missing Hodge-Riemann relations for $P_0^0$. (A
summand of a polarised Hodge structure is polarised.)

\subsection{Local Study of the Decomposition Theorem}

Our induction so far gives a decomposition
\begin{equation}
  \label{eq:lsdecomp}
  f_*\RM_X[n] \cong \bigoplus_{i \in \ZM} \pH^i(f_*\RM_X[n])[-i]
\end{equation}
and we know that each $\pH^i(f_*\RM_X[n])$ is semi-simple for $i \ne
0$. It remains to deduce the semi-simplicity of $\pH^0(f_*\RM_X[n])$. In
this section we outline the proof of the following:

\begin{prop} \label{prop:icsplitting}
There exist a local system $\LS_\l$ on each stratum $Y_\l$
  such that we have a canonical isomorphism
  \begin{equation}
\pH^0(f_*\RM_X[n]) = \bigoplus_{\l \in \Lambda} IC(\overline{Y_\lambda}, \LS_\l).\label{eq:icsplitting}
\end{equation}
\end{prop}

In the next section we explain why each $\LS_\l$ is semi-simple, which
completes the proof of the Semi-Simplicity Theorem, and hence of the
Decomposition Theorem.

As in the semi-small case one can reduce (by taking normal
slices) to the case of a point stratum $\{ y \}$. Denote by $i : \{ y
\} \into Y$ the inclusion. Let us say that
$\pH^0(f_*\RM_X[n])$ is \emph{semi-simple at $y$} if
\begin{equation} \label{eq:generalsplitting}
\pH^0(f_*\RM_X[n]) = i_*V \oplus \FS
\end{equation}
where $i_*V$ is the skyscraper sheaf at $y$ with stalk $V :=
H^0(\pH^0(f_*\RM_X[n])_y)$. (As in the semi-small case, $\FS$ is a
perverse sheaf whose structure can be ignored.)

\begin{rema}
  As in the semi-small case \eqref{eq:generalsplitting} is the key to
  establishing 
  \eqref{eq:icsplitting}. By induction we may assume that
  the restriction of $\FS$ to the complement of all point strata is a
  direct sum of intersection cohomology complexes as in
  \eqref{eq:icsplitting}, and \eqref{eq:generalsplitting} 
  allows us to extend this decomposition over point strata.
\end{rema}

Exactly as earlier we 
consider the form
\[
\Hom(i_*\RM_y, f_*\RM_X[n]) \times \Hom(f_*\RM_X[n],i_*\RM_y)
\to \End(i_*\RM_y) = \RM
\]
and again this form is canonically identified with the intersection
form
\begin{equation}
Q : H_n(f^{-1}(y)) \times H_n(f^{-1}(y)) \to \RM\label{eq:18}
\end{equation}
given by the embedding $f^{-1}(y) \into X$.

\begin{rema}
It will become clear in the discussion below 
that the rank of this form is precisely the multiplicity of $i_*\RM_y$
in $\pH^0(f_*\RM_X[n])$. However now the situation is considerably
more complicated because it is difficult to predict a priori what this
rank should be. (This should be contrasted with the semi-small case
where we knew that this multiplicity is always $\dim H_n(f^{-1}(y))$.)
As far as the Decomposition Theorem is
concerned, this is the fundamental difference between a
semi-small and a general map.  
\end{rema}


The technology of perverse sheaves provides a formal means of
circumventing this obstacle. Via the identifications (see Lemma
\ref{lem:caniso}) 
\begin{equation}
  \label{eq:17}
  H_n(f^{-1}(y)) = \Hom(i_*\RM_y, f_*\RM_X[n]) = H^0(i^!f_*\RM_X[n])
\end{equation}
the perverse filtration induces a filtration on $i^!f_*\RM_X[n]$, and
hence on $H_n(f^{-1}(y))$:
\begin{equation*}
  H_{n,\le i}(f^{-1}(y)) := \im ( H^0( i^! (\ptau_{\le i} f_* \RM_X[n])) \to H^0(i^!f_*\RM_X[n])).
\end{equation*}
We continue to refer to this as the \emph{perverse filtration}:
\begin{gather*}
 \dots\subset H_{n,\le -1}(f^{-1}(y)) \subset H_{n,\le 0}(f^{-1}(y))
\subset \dots
\end{gather*}

It turns out that the perverse filtration tells us precisely what the
radical of the intersection form should be.
From the definition of the perverse $t$-structure we have
$H^0(i^!(\ptau_{>0} f_*\RM_X[n])) = 0$. Hence:
\begin{equation}
  \label{eq:20}
  H_{n,\le 0}(f^{-1}(y)) = H_n(f^{-1}(y)).
\end{equation}
Also, as $\ptau_{<0} f_* \RM_X[n]$ does not contain any summand
isomorphic to $i_*\RM_y$ we deduce:
\begin{equation}
  \label{eq:21}
  H_{n, <0}(f^{-1}(y)) \subset \rad Q \subset H_n(f^{-1}(y)).
\end{equation}
Finally, $\pH^0(f^*\RM_X[n])$ is semi-simple at $y$ if and only the
composition 
\[
\Hom(i_*\RM_y, \pH^0(f_*\RM_X[n])) \times
\Hom(\pH^0(f_*\RM_X[n]) ,i_*\RM_y) \to \End(i_*\RM_y) = \RM
\]
is non-degenerate. We conclude:

\begin{prop}
$\pH^0(f_*\RM_X[n])$ is semi-simple at $y$ if and only if the
intersection form induces a non-degenerate form on
\[
H_{n,0}(f^{-1}(y)) := H_{n, \le 0}(f^{-1}(y)) / H_{n, <0}(f^{-1}(y)).
\]
\end{prop}

We use the $(\eta,\beta)$-Hodge-Riemann relations to conclude the proof:

\begin{theo}[Index Theorem for Maps]
The inclusion $f^{-1}(y) \into X$ yields an injection of pure Hodge structures
\[
H_{n,0}(f^{-1}(y)) \into H_0^0.
\]
This is an isometry with respect to the intersection form on the left
and the Poincar\'e form on the right. In particular, the intersection
form on $H_{n,0}(f^{-1}(y))$ underlies a polarization of pure Hodge
structure, and hence is non-degenerate.
\end{theo}

After taking the perverse filtrations into account, the proof (first
replace $Y$ by an affine neighbourhood of $y$, and then apply mixed
Hodge theory) is identical to the proof in the semi-small case.

\subsection{Semi-Simplicity of Local Systems}

It remains to see that all local systems occurring in the
decomposition
\begin{equation}
\pH^0(f_*\RM_X[n]) = \bigoplus_{\l \in \Lambda}
IC(\overline{Y}_\lambda, \LS_\l) \label{eq:decomp}
\end{equation}
are semi-simple. The idea is to exhibit each $\LS_\l$ (or more
precisely its restriction to a non-empty Zariski open subvariety $U
\subset Y_\l$) as the quotient of a local system associated with a
smooth proper map. Such local systems are semi-simple by Deligne's
theorem, and the semi-simplicity of each $\LS_\l$ follows. (Note that
$\pi_1(U) \onto \pi_1(Y_\l)$ is surjective, so it is enough to know
that the restriction of $\LS_\l$ to $U$ is semi-simple.)

Let $Y_\l \subset Y$ denote a stratum of dimension $s$. The local
system $\LS_\l$ occurring on the right of \eqref{eq:decomp} is
\[
\LS_\l := \mathcal{H}^{-s}(\pH^0(f_*\RM_X[n])_{Y_\l}).
\]
We must show that $\LS_\l$ is semi-simple.

By proper base
change
\[
H^j((f_*\RM_X[n])_y) = H^{n-j}(f^{-1}(y)) \quad \text{for all $y \in Y$.}
\]
The perverse filtration on $f_*\RM_X[n]$ induces a filtration on
$H(f^{-1}(y))$ which we denote by $H_{\le m}(f^{-1}(y))$. We set
\[
H_m(f^{-1}(y)) = H_{\le m}(f^{-1}(y)) / H_{< m}(f^{-1}(y)).
\]
For any $y \in Y_\l$ we have
\[
(\LS_\l)_y = H^{-s}( \pH^0(f_*\RM_X[n])_y) =  H^{n-s}_0(f^{-1}(y)).
\]
This equality exhibits the fibres of $\LS_\l$ as subquotients of
the cohomology of a variety. However we cannot apply Deligne's theorem
because the fibres $f^{-1}(y)$ are typically not smooth.

A key observation is that if $D \subset Y$ denotes the intersection
of $s$ general hyperplanes through $y$ which are transverse to all
strata then $f^{-1}(D)$ is smooth and
we have a surjection
\begin{equation} \label{eq:surj}
H^{n-s}(f^{-1}(D)) \onto H^{n-s}_0(f^{-1}(y)).
\end{equation}
This
observation is not difficult; the proof is similar to that of Lemma
\ref{lem:injective}.

With a little work (see \cite[\S 6.4]{dCM2}) one can find a
smooth family
\[
g : \XC_U \to U
\]
over a Zariski open subset $U \subset Y_\l$ whose
fibre over $y \in U$ is $f^{-1}(D)$, for some $D$ as in the previous
paragraph. The local system $R^{n-s} g_*\RM_{\XC_U}$ is semi-simple by
Deligne's theorem and one has a surjection 
$R^{n-s} g_*\RM_{\XC_U} \onto (\LS_\l)_{|U}$, which on stalks gives
maps as in \eqref{eq:surj}. The semi-simplicity of $\LS_\l$ follows.




\end{document}